\documentclass[english]{cccconf}
\usepackage[comma,numbers,square,sort&compress]{natbib}
\usepackage{epstopdf}
\usepackage{amsmath}
\usepackage{amsthm,diagbox}
\usepackage{amssymb}
\usepackage{subfigure}
\allowdisplaybreaks[4]

\begin{document}

\title{Distributed Nonsmooth Robust Resource Allocation with Cardinality Constrained Uncertainty}

\author{Yue Wei,
        Shuxin Ding, Hao Fang, Xianlin Zeng,
        Qingkai Yang, Bin Xin}

\affiliation{Beijing Institute of Technology,  Beijing 100081, P.~R.~China
        \email{yuewei@bit.edu.cn; jackietindd@gmail.com;  fangh@bit.edu.cn; xianlin.zeng@bit.edu.cn; qingkai.yang@bit.edu.cn; brucebin@bit.edu.cn}}
\maketitle

\begin{abstract}
A distributed nonsmooth robust resource allocation problem with cardinality constrained uncertainty is investigated in this paper. The global objective is consisted of local objectives, which are convex but nonsmooth. Each agent is constrained in its private convex set and has only the information of its corresponding local objective. The resource allocation condition is subject to the cardinality constrained uncertainty sets. By employing the duality theory of convex optimization, a dual problem of the robust resource allocation problem is presented. For solving this dual problem, a distributed primal-dual projected algorithm is proposed. Theoretically, the convergence analysis by using stability theory of differential inclusions is conducted. It shows that the algorithm can steer the multi-agent system to satisfy resource allocation condition at the optimal solution. In the end, a nontrivial simulation is shown and the results demonstrate the efficiency of the proposed algorithm.
\end{abstract}


\footnotetext{This work is supported by Projects of Major International (Regional) Joint Research Program NSFC (Grant no. 61720106011), NSFC (Grant no. 61621063, 61573062, 61673058, 61603378, 61603094).}

\section{Introduction}

In recent years, the distributed optimization problem is widely studied as a hot topic in the areas of machine learning \cite{Lasso} and multi-agent system coordination \cite{WC}-\cite{KSX}. In this problem, the objective is the sum of local objectives. Each agent can only obtain the knowledge of its private local objective. Many results of distributed optimization are focusing on steering the system to achieve consensus at the optimal solution \cite{C1}-\cite{C2}. On the other hand, the research of distributed globally constrained optimization has also gained a great of attention \cite{R1}-\cite{R2}, especially the distributed resource allocation problem. In order to solve the distributed resource allocation problem, 
a distributed gradient-based algorithm was proposed while the initialization of states is required \cite{IR}. After this work, the initialization-free distributed algorithms for distributed resource allocation have been investigated in \cite{RA1}-\cite{RA2}. 



While most of the existing works about distributed resource allocation have the assumption that the resource allocation condition is deterministic. This assumption may not applied for the distributed resource allocation problems applying in the real world. In order to solve these problems, robustness of the distributed resource allocation should be stressed. Robust optimization deals with uncertainty described by uncertain-but-bounded parameters \cite{RO}. Typically, there are several kinds of uncertain parameters (eg., box/interval uncertainty, ellipsoidal uncertainty, polyhedral uncertainty, cardinality constrained uncertainty, etc.) \cite{CCU}. Zeng et al. \cite{R1} proposed a distributed algorithm for robust resource allocation with polyhedral uncertain parameters. However, only considering polyhedral uncertain parameters may lead the problem too much conservative \cite{CON}. Cardinality constrained uncertainty provides a budget of uncertainty in terms of cardinality constraints which decrease the conservatism by combining interval and polyhedral uncertainty. Besides, many real-world robust optimization problems are related with cardinality constrained uncertainty \cite{CRW2}. Therefore, the robust optimization problem with cardinality constrained uncertainty needs to be analysed.



Nonsmooth optimization problem is increasingly popular due to its important role in a lot of signal processing, statistical inference and machine learning problems. In the compressed sensing problem, the sparsity-promoting regulator has the form of $l_{1}$-norm. In optimization problems with per-agent constraints, the indicator function of the constraint set of agent $i$ is nonsmooth. In the geometric median problem, the objective is the mean of a sum of $l_{2}$-norm functions.

In this paper, a distributed robust nonsmooth resource allocation problem with cardinality constrained uncertainty has been researched. The contributions of this paper are summed up as three parts:
\begin{enumerate}
	\item	The robust resource allocation problem we investigate here is with cardinality constrained uncertain parameters, which decrease the conservatism of the problem using polyhedral uncertain parameters.
	\item	We propose a distributed primal-dual projected algorithm with considering the duality theory of convex optimization.
	\item	The proof of the convergence of this algorithm has been given by employing the theory of nonsmooth analysis and differential inclusion.
\end{enumerate}
   
The paper is organized as follows. In Section 2, the necessary preliminary concepts of graph theory, projection operator and differential inclusion are introduced. Section 3 shows the robust nonsmooth resource allocation problem with cardinality constrained uncertainty. Section 4 proposes a distributed projected primal-dual algorithm. In Section 5 the convergence and correctness of the algorithm is proofed. Section 6 gives a numerical example to show the effectiveness of our proposed algorithm. Finally, Section 7 concludes this paper.

\section{Preliminary}
In this section, we introduce relevant notations, concepts on graph theory, projection operators and differential
inclusions. 

\subsection{Graph Theory}
A weighted undirected graph $\mathcal{G}$ is denoted by $\mathcal{G(V,E,A)}$, where $\mathcal{V} = \lbrace 1, \dots, n \rbrace$ is a set of nodes, $\mathcal{E} = \lbrace (i, k) : i, k \in \mathcal{V}; i \neq k \rbrace \subset \mathcal{V} \times \mathcal{V}$ is a set of edges, and $\mathcal{A} = [\alpha_{i,k}] \in \mathbb{R}^{n \times n}$ is a weighted adjacency matrix such that $\alpha_{i,k} = \alpha_{k,i} > 0$ if $(k,i) \in \mathcal{E}$ and $\alpha_{i,k} = 0$ otherwise, where $\mathbb{R}^{n \times n}$ denotes the set of $n$-by-$n$ real matrices. $k \in \mathcal{N}_{i}$ denotes agent $k$ is a neighbour of agent $i$. The Laplacian matrix is $L_{n} = D - \mathcal{A}$, where $D \in \mathbb{R}^{n \times n}$ is diagonal with $D_{i,i} = \sum^{n}_{k=1} \alpha_{i,k}$, $i \in \lbrace 1, \dots, n \rbrace$. Specifically, if the weighted graph $\mathcal{G}$ is undirected and connected, then $L_{n} = L^{T}_{n} \geq 0$.
\subsection{Projection Operator}

Define a projection operator as $P_{\Omega}(u) = \text{arg} \min_{v \in \Omega} \lbrace u-v \rbrace$, where $\Omega \subset \mathbb{R}^{n}$ is closed and convex, $\mathbb{R}^{n}$ denotes the set of $n$-dimensional real column vectors.

\textbf{Lemma 2.1.} \cite{PO} Let $\Omega \subset \mathbb{R}^{n}$ be closed and convex, and define $V : \mathbb{R}^{n} \to \mathbb{R}$ as $V(x) = \frac{1}{2} (\Vert x - P_{\Omega}(y)\Vert^{2} - \Vert x - P_{\Omega}(x)\Vert^{2})$ where $y \in \mathbb{R}^{n}$. Then $V(x) \geq \frac{1}{2} \Vert P_{\Omega}(x) - P_{\Omega}(y) \Vert^{2}$, $V(x)$ is differentiable and convex with respect to $x$, and $\nabla V(x) = P_{\Omega}(x) - P_{\Omega}(y)$.

\textbf{Lemma 2.2.} If $\Omega \subset \mathbb{R}^{n}$ is closed and convex, then $(P_{\Omega}(x) - P_{\Omega}(y))^{T}(x - y) \geq \Vert P_{\Omega}(x) - P_{\Omega}(y)\Vert^{2}$ for all $x, y \in \mathbb{R}^{n}$.

\subsection{Differential Inclusion}

Consider a nonsmooth system
\begin{equation}\label{Differential_Inclusion}
\dot{x} \in \mathcal{F}(x(t)), x(0) = x_{0}, t \geq 0
\end{equation}
where $\mathcal{F}: \mathbb{R}^{q} \to \mathcal{B}(\mathbb{R}^{q})$, $\mathcal{B}(\mathbb{R}^{q})$ is the collection of subsets of $\mathbb{R}^{q}$. A set $M$ is said to be weakly invariant (strongly invariant) with respect to (\ref{Differential_Inclusion})
if for any $x_{0} \in M$, $M$ contains a maximal solutions (all maximal solutions)) of (\ref{Differential_Inclusion}). An equilibrium point of (\ref{Differential_Inclusion}) is a point $x^{*} \in \mathbb{R}^{q}$ such that $\textbf{0}_{q} \in \mathcal{F}(x^{*})$. 

Let $V: \mathbb{R}^{q} \to \mathbb{R}$ be a locally Lipschitz continuous function and $\partial V$ be the Clarke generalized
gradient \cite{F} of $V(x)$ at $x$. The set-valued Lie derivative \cite{F} $\mathcal{L_{F}}V: \mathbb{R}^{q} \to \mathcal{B}(\mathbb{R})$ of $V$ with respect
to (\ref{Differential_Inclusion}) is defined as $\mathcal{L_{F}}V(x) = \lbrace a \in \mathbb{R}: \text{there exists } v \in \mathcal{F} \text{ such that } p^{T}v = a \text{ for all } p \in \partial V(x) \rbrace$. In the case when $\mathcal{L_{F}}V(x)$ is nonempty, we use $\max \mathcal{L_{F}}V(x)$ to denote the largest element of $\mathcal{L_{F}}V(x)$.

\textbf{Lemma 2.3.} \cite{CJ} For the differential inclusion ($\ref{Differential_Inclusion}$), we assume that $\mathcal{F}$ is upper semicontinuous and locally bounded, and $\mathcal{F}(x)$ takes nonempty, compact, and convex values. Let $V: \mathbb{R}_{q} \to \mathbb{R}$ be a locally Lipschitz and regular function, $S \subset \mathbb{R}_{q}$ be compact and strongly invariant for ($\ref{Differential_Inclusion}$),
$\phi(\cdot)$ be a solution of ($\ref{Differential_Inclusion}$), $\mathcal{R} = \lbrace x \in \mathcal{R}_{q}: 0 \in \mathcal{L_{F}}(x) \rbrace$, and $\mathcal{M}$ be the largest weakly invariant subset of $\bar{\mathcal{R}} \cap S$, where $\bar{\mathcal{R}}$ is the closure of $\mathcal{R}$. If $\max \mathcal{L_{F}}V(x) \leq 0$ for all $x \in S$, then $dist(\phi(t),\mathcal{M}) \to 0$ as $t \to \infty$.

\section{Problem Formulation}

In this section, the distribute nonsmooth robust resource allocation optimization problem is formulated. Consider the following distributed nonsmooth uncertain resource allocation problem 
\begin{align}
\min_{x_{i} \in \Omega_{i}}  & \quad f(x) = \sum_{i=1}^{n} f_{i}(x_{i}) \label{LO} \\
{\rm s.t.}  & \sum_{i=1}^{n} \bar{a}_{ij}^{l}x_{i}^{l}  \leq b_{j}^{l}, \quad \forall \bar{a}_{ij}^{l} \in \mathcal{U}_{ij}^{l} \notag
\end{align}
where $ j \in \lbrace 1, \cdots, m \rbrace, l \in \lbrace 1, \cdots, q \rbrace $, and with cardinality constrained uncertainty sets
\begin{gather}
\begin{split}
\mathcal{U}_{ij}^{l} = & \lbrace \bar{a}_{ij}^{l} \vert \bar{a}_{ij}^{l} \in [a_{ij}^{l} - \hat{a}_{ij}^{l}, a_{ij}^{l} +\hat{a}_{ij}^{l}], \\
& \sum\limits_{i,l}\!\left|  \frac{\bar{a}_{ij}^{l}-a_{ij}^{l}}{\hat{a}_{ij}^{l}}\right|  \leq \gamma_j, \forall i,j,l \rbrace
\end{split}
\end{gather}

For agent $i \in \lbrace 1, \cdots, n \rbrace$, $x_{i} \in \mathbb{R}^{q}$, $\Omega_{i}$ is the local constraint set, and $f_{i}(x_{i}) $ is the local objective which is continuous but not necessary smooth. $\bar{a}_{ij}^{l} \in \mathbb{R}$ is assumed to take arbitrary values in the uncertainty set $\mathcal{U}_{ij}^{l}$, $b_{j} = [b_{j}^{1}, \cdots, b_{j}^{q}]^{T} \in \mathbb{R}^{q}$, and $\gamma_{j}$ denotes the budget of uncertainty.

Then the corresponding robust optimization problem of the problem (\ref{LO}) is shown as
\begin{align}
\min_{x_{i} \in \Omega_{i}}  & \quad f(x) = \sum_{i=1}^{n} f^{i}(x_{i})  \notag \\
\text{s.t.}  & \sum_{i=1}^{n} a_{ij}^{l}x_{i}^{l} + \max_{S_{j}^{l} \in J_{j}^{l}: \vert S_{j}^{l} \vert = \gamma_{j}} \sum_{i \in S_{j}^{l}} \hat{a}_{ij}^{l} x_{i}^{l} \leq b_{j}^{l},  \notag \\
& j \in \lbrace 1, \cdots, m \rbrace, l \in \lbrace 1, \cdots, q \rbrace \label{Problem}
\end{align}

For the $l$-th dimensional elements of each agent's states with the $j$-th resource allocation condition, $S_{j}^{l}$ is a possible set of the chosen agents where the size of $S_{j}^{l}$ is $\gamma_{j}$, and $J_{j}^{l}$ is the set of all possible $S_{j}^{l}$.

According to the duality of convex optimization \cite{ROD}, the problem (\ref{Problem}) can be transferred to the corresponding dual problem as
\begin{gather}
\begin{split}
\min_{x_{i} \in \Omega_{i}}  & \quad f(x) = \sum_{i=1}^{n} f^{i}(x_{i})  \\
\text{s.t.} & \sum_{i=1}^{n} [A_{ij}x_{i} + \frac{1}{n} \gamma_{j} z_{ij} +  w_{ij} ] \leq \sum_{i=1}^{n} b_{ij},  \\
& \hat{A}_{ij} x_{i} \leq z_{ij} + w_{ij}, L_{mnq}Z = \textbf{0}_{mnq},  \\
& z_{ij} \geq \textbf{0}_{q}, w_{ij} \geq \textbf{0}_{q}, \\
& i \in \lbrace 1, \cdots, n \rbrace, j \in \lbrace 1, \cdots, m \rbrace
\end{split}\label{Dual-Problem}
\end{gather}
where $A_{ij} = diag \lbrace a_{ij}^{1}, \cdots, a_{ij}^{q}\rbrace \in \mathbb{R}^{q \times q}$, $\hat{A}_{ij} = diag \lbrace \hat{a}_{ij}^{1}, \cdots, \hat{a}_{ij}^{q}\rbrace \in \mathbb{R}^{q \times q}$, $z_{ij} \in \mathbb{R}^{q}$, $z_{i}=[(z_{i1})^{T},\dots,(z_{im})^{T}]^{T}$, $Z=[(z_{1})^{T},\dots,(z_{n})^{T}]^{T}$, $w_{ij} \in \mathbb{R}^{q}$. $\sum_{i=1}^{n} b_{ij} = b_{j}$, $L_{mnq} = L_{n} \otimes I_{mq} $, where $A\otimes B$ denotes the the Kronecker product of matrices $A$ and $B$. $\textbf{0}_{n}$ is the $n \times 1$ vector with all elements of 0.

The assumptions below are made for the wellposedness of the problem ($\ref{Dual-Problem}$) in this section.

\textbf{Assumption 3.1.}  1) The weighted graph $\mathcal{G}$ is connected and undirected.

2) For $i \in \lbrace 1, \cdots, n\rbrace$, $f_{i}$ is strictly convex on an open set containing $\Omega_{i}$, and $\Omega_{i} \subset \mathbb{R}^{q}$ is closed and convex.

3) (Slater's constraint condition) There exist $x_{i} \in \Omega_{i}$, $z_{ij} \in \bar{\mathbb{R}}_{+}^{q}$ and $w_{ij} \in \bar{\mathbb{R}}_{+}^{q}$ satisfying the constraint for $i \in \lbrace 1, \cdots, n \rbrace$ and $j \in \lbrace 1, \cdots, m \rbrace$, where $\bar{\mathbb{R}}_{+}^{q}$ denotes the set of nonnegative $q$-dimensional real column vectors.

The Lagrangian of dual problem (\ref{Dual-Problem}) is described as
\begin{align}
&L(x,Z,W,\Lambda^{1},\Lambda^{2}, U) \notag \\
= & \sum_{i=1}^{n} [f^{i}(x_{i}) + \sum_{j=1}^{m} [(\lambda^{1}_{ij})^{T} (A_{ij}x_{i} + \frac{1}{n}\gamma_{j} z_{ij} + w_{ij} - b_{ij}) \notag \\
&+ (\lambda^{2}_{ij})^{T}(\hat{A}_{ij}x_{i} - z_{ij} - w_{ij}) ]] + \mu^{T} L_{mnq}Z
\end{align}
where $w_{i}=[(w_{i1})^{T},\dots,(w_{im})^{T}]^{T}$, 
$W=[(w_{1})^{T},\dots,(w_{n})^{T}]^{T}$, 
$\lambda_{i}^{g}=[(\lambda_{i1}^{g})^{T},\dots,(\lambda_{im}^{g})^{T}]^{T}$, $\Lambda^{g}=[(\lambda_{1}^{g})^{T},\dots,(\lambda_{n}^{g})^{T}]^{T}$, $\mu_{i}=[(\mu_{i1})^{T},\dots,(\mu_{im})^{T}]^{T}$,  $U=[\mu_{1}^{T},\dots,\mu_{n}^{T}]^{T}$, $g \in \lbrace 1,2 \rbrace$.

Then according to problem (\ref{Dual-Problem}), the following lemma is arrived by the Karush-Kuhn-Tucker (KKT) condition of convex optimization problems.

\textbf{Lemma 3.1.}
Under the Assumptions 3.1, a feasible point $x^{*} \in \mathbb{R}^{nq}$ is a minimizer to Problem (\ref{Dual-Problem}) if and only if there exist $x^{*}_{i} \in \Omega_{i} \in \mathbb{R}^{q}$, $\lambda^{1*}_{ij} \in \mathbb{R}^{q}$, $\lambda^{2*}_{ij} \in \mathbb{R}^{q}$, $\mu_{ij}^{*} \in \mathbb{R}^{q}$ , $w_{ij}^{*} \in \mathbb{R}^{q}$ and $z_{ij}^{*} \in \mathbb{R}^{q}$ such that for $i \in \lbrace 1, \cdots, n \rbrace$, $j \in \lbrace 1, \cdots, m \rbrace$
{\setlength\abovedisplayskip{-3pt}
\setlength\belowdisplayskip{1pt}
\begin{subequations}
\begin{align}
-\partial f_{i}(x_{i}^{*}) - \sum_{j=1}^{m} A_{ij} \lambda^{1*}_{ij}  - \sum_{j=1}^{m}  \hat{A}_{ij} \lambda^{2*}_{ij} \in & \mathcal{N}_{\Omega_{i}}(x^{*}_{i}), \label{K1}\\
- \frac{1}{n} \gamma_{j} \lambda^{1*}_{ij} + \lambda^{2*}_{ij} - \sum_{k\in \mathcal{N}_{i}} (\mu_{ij} - \mu_{kj}) \in & \mathcal{N}_{\bar{\mathbb{R}}_{+}^{q}}(z^{*}_{ij}),  \label{K2}\\
- \lambda^{1*}_{ij} + \lambda^{2*}_{ij} \in & \mathcal{N}_{\bar{\mathbb{R}}_{+}^{q}}(w^{*}_{ij}),  \label{K3}\\
\sum_{i=1}^{n} [A_{ij}x_{i}^{*} + \frac{1}{n}\gamma_{j} z_{ij}^{*} + w_{ij}^{*} - b_{ij}] \leq & \textbf{0}_{q}, \label{K4}\\
\sum_{i=1}^{n} [\hat{A}_{ij}x_{i}^{*} - z_{ij}^{*} - w_{ij}^{*}] \leq & \textbf{0}_{q}, \label{K5}\\
L_{mnq}Z^{*} = & \textbf{0}_{mnq}, \label{K6}\\
(\lambda_{ij}^{1*})^{T} [A_{ij}x_{i}^{*} + \frac{1}{n}\gamma_{j} z_{ij}^{*} + w_{ij}^{*} - b_{ij}] = & 0, \label{K7}\\
(\lambda_{ij}^{2*})^{T} [\hat{A}_{ij}x_{i}^{*} - z_{ij}^{*} - w_{ij}^{*}] = & 0. \label{K8}
\end{align}\label{KKT}%
\end{subequations}
}%
where $\mathcal{N}_{\Omega_{i}}(x^{*}_{i})$ is the the normal cone of $\Omega_{i}$ at $x^{*}_{i}$.

The proof of Lemma 3.1 is omitted since it is a trivial extension of the proof for Theorem 3.34 in \cite{NO}. 

\section{Algorithm Design}
In this section, we propose a distributed algorithm for this problem (\ref{Dual-Problem}). The algorithm is detailed as below:

\begin{gather}
\begin{cases}
\dot{\bar{x}}_{i} \in & \!\!\!\!\!\! - \bar{x}_{i} \!\!+\!\! x_{i} \!\!-\!\! \partial f_{i}(x_{i}) \!\!-\!\! \sum_{j=1}^{m} \!\! A_{ij} \lambda^{1}_{ij} \!\!-\!\! \sum_{j=1}^{m} \!\! \hat{A}_{ij} \lambda^{2}_{ij} \\
\dot{\bar{z}}_{ij} = & \!\!\!\!\!\!- \bar{z}_{ij} \!\!+\!\! z_{ij} \!\!-\!\! \frac{1}{n} \gamma_{j} \lambda^{1}_{ij} \!\!+\!\! \lambda^{2}_{ij} \!\!-\!\! \sum_{k\in \mathcal{N}_{i}} \alpha_{k}(\mu_{ij} \!\!-\!\! \mu_{kj}) \\
\dot{\bar{w}}_{ij} = & \!\!\!\!\!\!- \bar{w}_{ij} + w_{ij} - \lambda^{1}_{ij} + \lambda^{2}_{ij} \\
\dot{\mu}_{ij} = & \!\!\!\!\!\!\sum_{k\in \mathcal{N}_{i}}\alpha_{ik}(z_{ij} - z_{kj}) \\
\dot{\bar{\lambda}}_{ij}^{1} = & \!\!\!\!\!\!- \bar{\lambda}_{ij}^{1} + \lambda_{ij}^{1} + [A_{ij}x_{i} + \frac{1}{n}\gamma_{j} z_{ij} + w_{ij} - b_{ij}] \\
& \!\!\!\!\!\!+\!\! \sum_{k\in \mathcal{N}_{i}}\!\!\alpha_{ik}(y_{ij}^{1} \!\!-\!\! y_{kj}^{1}) \!\!-\!\! \sum_{k\in \mathcal{N}_{i}}\!\!\alpha_{ik}(\lambda_{ij}^{1} \!\!-\!\! \lambda_{kj}^{1}) \\
\dot{\bar{\lambda}}_{ij}^{2} = & \!\!\!\!\!\!-\bar{\lambda}_{ij}^{2} + \lambda_{ij}^{2} + [\hat{A}_{ij}x_{i} - z_{ij} - w_{ij}] \\
& \!\!\!\!\!\!+\!\! \sum_{k\in \mathcal{N}_{i}}\!\!\alpha_{ik}(y_{ij}^{2} \!\!-\!\! y_{kj}^{2}) \!\!-\!\! \sum_{k\in \mathcal{N}_{i}}\!\!\alpha_{ik}(\lambda_{ij}^{2} \!\!-\!\! \lambda_{kj}^{2}) \\
\dot{y}_{ij}^{1} = &\!\!\!\!\!\! - \sum_{k\in \mathcal{N}_{i}}\alpha_{ik}(\lambda_{ij}^{1} - \lambda_{kj}^{1}) \\
\dot{y}_{ij}^{2} = &\!\!\!\!\!\! - \sum_{k\in \mathcal{N}_{i}}\alpha_{ik}(\lambda_{ij}^{2} - \lambda_{kj}^{2}) \\
x_{i} = &\!\!\!\!\!\! P_{\Omega_{i}}[\bar{x}_{i}], z_{ij} = P_{\bar{\mathbb{R}}_{+}^{q}}[\bar{z}_{ij}], w_{ij} = P_{\bar{\mathbb{R}}_{+}^{q}}[\bar{w}_{ij}], \\
\lambda_{ij}^{1} = &\!\!\!\!\!\! P_{\bar{\mathbb{R}}_{+}^{q}}[\bar{\lambda}_{ij}^{1}], \lambda_{ij}^{2} = P_{\bar{\mathbb{R}}_{+}^{q}}[\bar{\lambda}_{ij}^{2}] 
\end{cases} \label{Algorithm}
\end{gather}
where $t \geq 0$.

The algorithm (\ref{Algorithm}) can be also written as a compact form as
\begin{gather}
\begin{split}
\dot{\Phi} \in & \mathcal{F}(\Phi), x = P_{\Omega}[\bar{x}], Z = P_{\bar{\mathbb{R}}_{+}^{mnq}}[\bar{Z}], \\
W = & P_{\bar{\mathbb{R}}_{+}^{mnq}}[\bar{W}], \Lambda^{1} = P_{\bar{\mathbb{R}}_{+}^{mnq}}[\bar{\Lambda}^{1}], \Lambda^{2} = P_{\bar{\mathbb{R}}_{+}^{mnq}}[\bar{\Lambda}^{2}]
\end{split}\label{CAlgorithm}
\end{gather}
where $\Phi = [\bar{x}^{T}, \bar{Z}^{T}, \bar{W}^{T}, U^{T}, (\bar{\Lambda}^{1})^{T}, (\bar{\Lambda}^{2})^{T}, (Y^{1})^{T}, (Y^{2})^{T}]^{T}$, 
$P_{\Omega}[\bar{x}]=[(P_{\Omega_{1}}[\bar{x}_{1}])^{T}, \cdots, (P_{\Omega_{n}}[\bar{x}_{n}])^{T}]^{T}$,
$\bar{z}_{i}=[(\bar{z}_{i1})^{T},\dots,(\bar{z}_{im})^{T}]^{T}$, 
$\bar{Z}=[(\bar{z}_{1})^{T},\dots,(\bar{z}_{n})^{T}]^{T}$, 
$\bar{w}_{i}=[(\bar{w}_{i1})^{T},\dots,(\bar{w}_{im})^{T}]^{T}$, 
$\bar{W}=[(\bar{w}_{1})^{T},\dots,(\bar{w}_{n})^{T}]^{T}$, 
$\bar{\lambda}_{i}^{g}=[(\bar{\lambda}_{i1}^{g})^{T},\dots,(\bar{\lambda}_{im}^{g})^{T}]^{T}$, $\bar{\Lambda}^{g}=[(\bar{\lambda}_{1}^{g})^{T},\dots,(\bar{\lambda}_{n}^{g})^{T}]^{T}$, 
$y_{i}^{g}=[(y_{i1}^{g})^{T},\dots,(y_{im}^{g})^{T}]^{T}$, 
$Y^{g} = [(y_{1}^{g})^{T},\dots,(y_{n}^{g})^{T}]^{T}$, $g \in \lbrace 1,2 \rbrace$.

In (\ref{CAlgorithm}), $\mathcal{F}(\phi)$ is defined as 
\begin{gather}
\begin{split}
\mathcal{F}(\phi) = & \lbrace [p_{\bar{x}}^{T}, p_{\bar{Z}}^{T}, p_{\bar{W}}^{T}, p_{U}^{T}, p_{\bar{\Lambda}^{1}}^{T}, p_{\bar{\Lambda}^{2}}^{T}, p_{Y^{1}}^{T}, p_{ Y^{2}}^{T}]^{T}   \\
\in &  \mathbb{R}^{nq} \times \mathbb{R}^{mnq} \times \mathbb{R}^{mnq} \times \mathbb{R}^{mnq}  \\
& \times \mathbb{R}^{mnq} \times \mathbb{R}^{mnq} \times \mathbb{R}^{mnq} \times \mathbb{R}^{mnq} \rbrace  
\end{split}
\end{gather}
with
\begin{gather}
\begin{cases}
p_{\bar{x}} = & \!\!\!\!\!\!- \bar{x} + x - f_{x} - EA^{*}\Lambda^{1} - E\hat{A}^{*}\Lambda^{2}  \\
p_{\bar{Z}} = & \!\!\!\!\!\!- \bar{Z} + Z - \frac{1}{n}\Gamma \Lambda^{1} + \Lambda^{2} - L_{mnq} U  \\
p_{\bar{W}} = & \!\!\!\!\!\!- \bar{W} + W - \Lambda^{1} + \Lambda^{2}, p_{U} = L_{mnq}Z  \\
p_{\bar{\Lambda}^{1}} = & \!\!\!\!\!\!- \bar{\Lambda}^{1} + \Lambda^{1} + [A^{*}E^{T}x + \frac{1}{n}\Gamma Z + W - B]  \\
& \!\!\!\!\!\!+ L_{mnq}Y^{1} - L_{mnq}\Lambda^{1}  \\
p_{\bar{\Lambda}^{2}} = & - \bar{\Lambda}^{2} + \Lambda^{2} + [\hat{A}^{*}E^{T}x - Z - W]  \\
& \!\!\!\!\!\!+ L_{mnq}Y^{2} - L_{mnq}\Lambda^{2}  \\
p_{Y^{1}} = & - L_{mnq}\Lambda^{1}, p_{Y^{2}} = - L_{mnq}\Lambda^{2} 
\end{cases}
\end{gather}
where $E = I_{n} \otimes (\textbf{1}_{m} \otimes I_{q}) \in \mathbb{R}^{nq \times mnq}$, $I_{n}$ is the $n$-dimensional identity matrix, $\textbf{1}_{m}$ denotes the $n \times 1$ vector with all elements of 1. $A_{i} = diag\lbrace A_{i1}, \cdots, A_{im}\rbrace \in \mathbb{R}^{mq \times mq}$, $A^{*} = diag\lbrace A_{1}, \cdots, A_{n}\rbrace \in \mathbb{R}^{mnq \times mnq}$, $\hat{A}_{i} = diag\lbrace \hat{A}_{i1}, \cdots, \hat{A}_{im}\rbrace \in \mathbb{R}^{mq \times mq}$, $\hat{A}^{*} = diag\lbrace \hat{A}_{1}, \cdots, \hat{A}_{n}\rbrace \in \mathbb{R}^{mnq \times mnq}$, $\Gamma = I_{n} \otimes diag\lbrace \gamma_{1}I_{q}, \cdots, \gamma_{m}I_{q}\rbrace \in \mathbb{R}^{mnq \times mnq}$, $L_{mnq} = L_{n} \otimes I_{mq}$, $b_{i}=[(b_{i1})^{T},\dots,(b_{im})^{T}]^{T}$, $B = [(b_{1})^{T},\dots,(b_{n})^{T}]^{T}$, $f_{x} \in \partial f(x)$.

Then the equilibrium of algorithm (\ref{CAlgorithm}) is
\begin{subequations}
\begin{align}
\textbf{0}_{nq} \!\!= & - \bar{x}^{*} + x^{*} - f_{x^{*}} - EA^{*}\Lambda^{1*} - E\hat{A}^{*}\Lambda^{2*} \label{E1} \\
\textbf{0}_{mnq} \!\!= & - \bar{Z}^{*} + Z^{*} - \frac{1}{n}\Gamma \Lambda^{1*} + \Lambda^{2*} - L_{mnq} U^{*}  \label{E2}\\
\textbf{0}_{mnq} \!\!= & - \bar{W}^{*} + W^{*} - \Lambda^{1*} + \Lambda^{2*} \label{E3}\\
\textbf{0}_{mnq} \!\!= & L_{mnq}Z^{*}  \label{E4}\\
\textbf{0}_{mnq} \!\!= & \!\!-\!\! \bar{\Lambda}^{1*} \!\!\!\!+\!\! \Lambda^{1*} \!\!\!\!+\!\! [\!A^{*}\!E^{T}\!\!x^{*}\!\!\!\! +\!\! \frac{1}{n}\Gamma\! Z^{*}\!\!\!\! +\!\! W^{*}\!\!\! -\!\! B] \!\!+\!\! L_{mnq}Y^{1*}  \label{E5} \\
\textbf{0}_{mnq} \!\!= & \!-\! \bar{\Lambda}^{2*} \!\!+\! \Lambda^{2*} \!\!+\! [\hat{A}^{*}\!E^{T}\!\!x^{*} \!\!\!-\!\! Z^{*} \!\!\!-\! W^{*}]\!\!+\! L_{mnq}Y^{2*}  \label{E6} \\
\textbf{0}_{mnq} \!\!= & L_{mnq}\Lambda^{1*}, \textbf{0}_{mnq} = L_{mnq}\Lambda^{2*}  \label{E7} \\
x^{*} \!\!= & P_{\Omega}[\bar{x}^{*}], Z^{*} \!=\! P_{\bar{\mathbb{R}}_{+}^{mnq}}[\bar{Z}^{*}], W^{*}\! =\! P_{\bar{\mathbb{R}}_{+}^{mnq}}[\bar{W}^{*}]  \label{E8}\\
\Lambda^{1*} \!\!= & P_{\bar{\mathbb{R}}_{+}^{mnq}}[\bar{\Lambda}^{1*}], \Lambda^{2*} = P_{\bar{\mathbb{R}}_{+}^{mnq}}[\bar{\Lambda}^{2*}] \label{E9}
\end{align}\label{Equilibrium}
\end{subequations}

Here we give the Lemma 4.1 to link the equilibrium of algorithm with the solution of problem (\ref{Problem}).

\textbf{Lemma 4.1.}
Consider Problem ($\ref{Problem}$) and Assumption 3.1 holds. If $ \phi^{*} \in \mathbb{R}^{(7m+1)nq} $ is an equilibrium of ($\ref{Algorithm}$), then $x^{*} = P_{\Omega}[\bar{x}^{*}]$ is a solution to Problem ($\ref{Problem}$). 
\begin{proof}
Suppose $\phi^{*} \in \mathbb{R}^{(7m+1)nq}$ is an equilibrium of ($\ref{Algorithm}$).  When considering ($\ref{E1}$), ($\ref{E2}$) and ($\ref{E3}$), there exists $f_{x^{*}} \in \partial f(x^{*})$ such that $\bar{x}^{*} = x^{*} - f_{x^{*}} - EA^{*}\Lambda^{1*} - E\hat{A}^{*}\Lambda^{2*}$, $\bar{Z}^{*} = Z^{*} - \frac{1}{n}\Gamma \Lambda^{1*} + \Lambda^{2*} - L_{mnq} U^{*}$, $\bar{W}^{*} = W^{*} - \Lambda^{1*} + \Lambda^{2*}$. Since $x^{*} = P_{\Omega}[\bar{x}^{*}]$, $Z^{*} = P_{\bar{\mathbb{R}}_{+}^{mnq}}[\bar{Z}^{*}]$, $W^{*} = P_{\bar{\mathbb{R}}_{+}^{mnq}}[\bar{W}^{*}]$, it follows that ($\ref{K1}$), ($\ref{K2}$) and ($\ref{K3}$) holds.

According to (\ref{E5}) and (\ref{E6}), one can have that
\begin{align}
&Q_{j} (- \bar{\Lambda}^{1*} + \Lambda^{1*} + [A^{*}E^{T}x^{*} + \frac{1}{n}\Gamma Z^{*} + W^{*} - B] \notag\\
&+ L_{mnq}Y^{1*}) \!=\! - \sum_{i=1}^{n} (\bar{\lambda}^{1*}_{ij} \!-\! \lambda^{1*}_{ij}) \!+\! \sum_{i=1}^{n} H_{ij}^{1*} = \textbf{0}_{q} \label{H1}\\  
&Q_{j} (- \bar{\Lambda}^{2*} + \Lambda^{2*} + [\hat{A}^{*}E^{T}x^{*} - Z^{*} - W^{*}] + L_{mnq}Y^{2*}) \notag \\
&= - \sum_{i=1}^{n} (\bar{\lambda}^{2*}_{ij} - \lambda^{2*}_{ij}) + \sum_{i=1}^{n} H_{ij}^{2*} = \textbf{0}_{q}  \label{H2}
\end{align}
where $Q_{j} = \textbf{1}_{n}^{T} \otimes (I_{m}^{j} \otimes I_{q}) \in \mathcal{R}^{q \times mnq}$, $I_{m}^{j}$ denotes the $j$-th row of $I_{m}$. $Q_{j} L_{mnq} Y^{1*} = \textbf{0}_{q}$, $Q_{j} L_{mnq} Y^{2*} = \textbf{0}_{q}$, $H^{1}_{i,j} = A_{ij}x_{i} + \frac{1}{n}\gamma_{j} z_{ij} + w_{ij} - b_{ij}$, $H^{2}_{i,j} = \hat{A}_{ij}x_{i} - z_{ij} - w_{ij}$. Since $\lambda_{ij}^{1} = P_{\bar{\mathbb{R}}_{+}^{mnq}}[\bar{\lambda}_{ij}^{1}] \geq \textbf{0}_{q}$, $\lambda_{ij}^{2} = P_{\bar{\mathbb{R}}_{+}^{mnq}}[\bar{\lambda}_{ij}^{2}] \geq \textbf{0}_{q}$, $\bar{\lambda}_{ij}^{1} - \lambda_{ij}^{1} \leq \textbf{0}_{q}$ and $\bar{\lambda}_{ij}^{2} - \lambda_{ij}^{2} \leq \textbf{0}_{q}$ for all $i \in \lbrace 1, \cdots, n \rbrace$. Hence ($\ref{K4}$) and ($\ref{K5}$) holds. ($\ref{E4}$) equals to ($\ref{K6}$), which means that ($\ref{K6}$) holds.


It follows from ($\ref{E7}$), ($\ref{E9}$) and $\lambda_{ij}^{1} = P_{\bar{\mathbb{R}}_{+}^{q}}[\bar{\lambda}_{ij}^{1}] \geq \textbf{0}_{q}$, $\lambda_{ij}^{2} = P_{\bar{\mathbb{R}}_{+}^{q}}[\bar{\lambda}_{ij}^{2}] \geq \textbf{0}_{q}$ that there exist $\lambda^{1*}_{0} \in \bar{\mathbb{R}}_{+}^{mq}$ and $\lambda^{2*}_{0} \in \bar{\mathbb{R}}_{+}^{mq}$ such that $\Lambda^{1*} = \lambda^{1*}_{0} \otimes \textbf{1}_{n}$ and $\Lambda^{2*} = \lambda^{2*}_{0} \otimes \textbf{1}_{n}$. If $\lambda^{1*}_{0} = \textbf{0}_{mq}$ and $\lambda^{2*}_{0} = \textbf{0}_{mq}$, the ($\ref{K7}$) and ($\ref{K8}$) holds. If $\lambda^{1*}_{0} > \textbf{0}_{mq}$ and $\lambda^{2*}_{0} > \textbf{0}_{mq}$, it is clear that $\hat{\lambda}_{ij}^{1*} = \lambda_{ij}^{1*}$, $\hat{\lambda}_{ij}^{2*} = \lambda_{ij}^{2*}$, $\sum_{i=1}^{n} H^{1*}_{ij} = \textbf{0}_{q}$ and $\sum_{i=1}^{n} H^{2*}_{ij} = \textbf{0}_{q}$, which means that ($\ref{K7}$) and ($\ref{K8}$) also holds.

By Lemma 3.1, $(x^{*}, Z^{*}, W^{*}) \in \Omega \times \mathbb{R}^{mnq} \times \mathbb{R}^{mnq}$ is an optimal solution of Problem ($\ref{Dual-Problem}$). Note that Problem ($\ref{Dual-Problem}$) is the strong dual problem of Problem ($\ref{Problem}$). Then the proof is accomplished.
\end{proof}

\section{Main Result}
In this section, we give the convergence analysis of 
our algorithm (\ref{Algorithm}). Define the Lyapunov candidate
\begin{gather}
\begin{split}
V(\phi) = & V_{1}(\bar{x}) + V_{2}(\bar{Z}) + V_{3}(\bar{W}) + V_{4}(U) \\
& + V_{5}(\bar{\Lambda}^{1}) + V_{6}(\bar{\Lambda}^{2}) + V_{7}(Y^{1}) + V_{8}(Y^{2}) 
\end{split} \label{Lyapunov}
\end{gather}
where
\begin{gather}
\begin{cases}
V_{1}(\bar{x}) = \frac{1}{2} (\Vert \bar{x} - x^{*} \Vert^{2} - \Vert \bar{x} - x \Vert^{2}) \\
V_{2}(\bar{Z}) = \frac{1}{2} (\Vert \bar{Z} - Z^{*} \Vert^{2} - \Vert \bar{Z} - Z \Vert^{2}) \\
V_{3}(\bar{W}) = \frac{1}{2} (\Vert \bar{W} - W^{*} \Vert^{2} - \Vert \bar{W} - W \Vert^{2}) \\
V_{4}(U) = \frac{1}{2} (\Vert U - U^{*} \Vert^{2}) \\
V_{5}(\bar{\Lambda}^{1}) = \frac{1}{2} (\Vert \bar{\Lambda}^{1} - \Lambda^{1*} \Vert^{2} - \Vert (\Vert \bar{\Lambda}^{1} - \Lambda^{1} \Vert^{2}) \\
V_{6}(\bar{\Lambda}^{2}) = \frac{1}{2} (\Vert \bar{\Lambda}^{2} - \Lambda^{2*} \Vert^{2} - \Vert (\Vert \bar{\Lambda}^{2} - \Lambda^{2} \Vert^{2}), \\
V_{7}(Y^{1}) = \frac{1}{2} (\Vert Y^{1} - Y^{1*} \Vert^{2}) \\
V_{8}(Y^{2}) = \frac{1}{2} (\Vert Y^{2} - Y^{2*} \Vert^{2})
\end{cases}
\end{gather}

In the following lemma, we have analysed the set-valued derivative of $V(\phi)$ defined in ($\ref{Lyapunov}$) along the trajectories of Algorithm ($\ref{Algorithm}$).

\textbf{Lemma 5.1.} Consider Algorithm ($\ref{Algorithm}$) under Assumption 3.1 with $V(\phi)$ defined in ($\ref{Lyapunov}$). If $\beta \in \mathcal{L_{F}} V(\phi)$, then there exist $f_{x} \in \partial f(x)$ and $f_{x^{*}} \in \partial f(x^{*})$ with $x = P_{\Omega}[\bar{x}]$ and $x^{*} = P_{\Omega}[\bar{x}^{*}]$ such that
\begin{gather}
\begin{split}
\beta \leq & - (x - x^{*})^{T}(f_{x} - f_{x^{*}}) - (\Lambda^{1})^{T}L_{mnq}\Lambda^{1}  \\
& - (\Lambda^{2})^{T}L_{mnq}\Lambda^{2} \leq 0
\end{split} \label{BA}
\end{gather}

\begin{proof}
It follows from Lemma 2.1 that the gradients of $V(\phi)$ with respect to $\phi$ are 
\begin{gather}
\begin{cases}
\nabla_{\bar{x}} V(\phi) = x - x^{*},\nabla_{\bar{Z}} V(\phi) = Z - Z^{*}\\
\nabla_{\bar{W}} V(\phi) = W - W^{*},\nabla_{U} V(\phi) = U - U^{*}\\
\nabla_{\bar{\Lambda}^{1}} V(\phi) = \Lambda^{1} - \Lambda^{1*}, \nabla_{\bar{\Lambda}^{2}} V(\phi) = \Lambda^{2} - \Lambda^{2*}\\
\nabla_{Y^{1}} V(\phi) = Y^{1} \! - \! Y^{1*}, \nabla_{Y^{2}} V(\phi) = Y^{2} - Y^{2*}
\end{cases} 
\end{gather}

The function $V(\phi)$ with the trajectories of ($\ref{Algorithm}$) satisfies

\begin{gather}
\begin{split}
\mathcal{L_{F}} V(\phi) = & \lbrace \!\beta \!\in\! \mathbb{R} \!:\! \beta \!=\! \nabla_{\bar{x}} V(\phi)^{T}p_{\bar{x}} \!+\! \nabla_{\bar{Z}} V(\phi)^{T}p_{\bar{Z}} \\
& + \nabla_{\bar{x}} V(\phi)^{T}p_{\bar{W}} + \nabla_{U} V(\phi)^{T}p_{U}  \\
& + \nabla_{\bar{\Lambda}^{1}} V(\phi)^{T}p_{\bar{\Lambda}^{1}} + \nabla_{\bar{\Lambda}^{2}} V(\phi)^{T}p_{\bar{\Lambda}^{2}}  \\
& + \nabla_{Y^{1}} V(\phi)^{T}p_{Y^{1}} + \nabla_{Y^{2}} V(\phi)^{T}p_{Y^{2}} \rbrace
\end{split} 
\end{gather}

Suppose $\beta \in \mathcal{L_{F}}V(\phi)$. There exists $f_{x} \in \partial f(x)$ such that $\beta = \sum_{i=1}^{8} \beta_{i}$, then 

\begin{gather}\label{Beta_8}
\begin{cases}
\beta_{1} = &\!\!\!\!\!\! (x - x^{*})^{T}(- \bar{x} + x \!-\! f_{x} \!\!-\!\! EA^{*}\Lambda^{1} \!\!-\!\! E\hat{A}^{*}\Lambda^{2})
\\
\beta_{2} = &\!\!\!\!\!\! (Z \!-\! Z^{*})^{T}(- \bar{Z} \!+\! Z \!-\! \frac{1}{n}\Gamma \Lambda^{1} \!+\! \Lambda^{2} \!-\! L_{mnq} U)\\
\beta_{3} = &\!\!\!\!\!\! (W - W^{*})^{T}(- \bar{W} + W - \Lambda^{1} + \Lambda^{2})\\
\beta_{4} = &\!\!\!\!\!\! (U - U^{*})^{T}L_{mnq}Z\\
\beta_{5} = &\!\!\!\!\!\! (\Lambda^{1} - \Lambda^{1*})^{T}(- \bar{\Lambda}^{1} + \Lambda^{1} + [A^{*}E^{T}x + \frac{1}{n}\Gamma Z \\
&\!\!\!\!\!\! + W - B] + L_{mnq}Y^{1} - L_{mnq}\Lambda^{1})
\\
\beta_{6} = &\!\!\!\!\!\! (\Lambda^{2} - \Lambda^{2*})^{T}(- \bar{\Lambda}^{2} + \Lambda^{2} + [\hat{A}^{*}E^{T}x - Z \\
&\!\!\!\!\!\! - W] + L_{mnq}Y^{2} - L_{mnq}\Lambda^{2})
\\
\beta_{7} = &\!\!\!\!\!\! -(Y^{1} - Y^{1*})^{T}L_{mnq}\Lambda^{1}
\\
\beta_{8} = & -(Y^{2} - Y^{2*})^{T}L_{mnq}\Lambda^{2}
\end{cases}
\end{gather}

Since $ \phi^{*} \in \mathbb{R}^{(7m+1)nq}$ is an equilibrium of ($\ref{Algorithm}$), there exists $f_{x^{*}} \in \partial f(x^{*})$ such that ($\ref{Equilibrium}$) holds.

From ($\ref{Equilibrium}$) and ($\ref{Beta_8}$), one can have that
\begin{gather}
\begin{split}
\beta = & - (x - x^{*})^{T}(\bar{x} - \bar{x}^{*}) + \Vert x - x^{*} \Vert^{2} \\
& - (Z - Z^{*})^{T}(\bar{Z} - \bar{Z}^{*}) + \Vert Z - Z^{*} \Vert^{2} \\
& - (W - W^{*})^{T}(\bar{W} - \bar{W}^{*}) + \Vert W - W^{*} \Vert^{2} \\
& - (\Lambda^{1} - \Lambda^{1*})^{T}(\bar{\Lambda}^{1} - \bar{\Lambda}^{1*}) + \Vert \Lambda^{1} - \Lambda^{1*} \Vert^{2} \\
&  - (\Lambda^{2} - \Lambda^{2*})^{T}(\bar{\Lambda}^{2} - \bar{\Lambda}^{2*}) + \Vert \Lambda^{2} - \Lambda^{2*} \Vert^{2} \\
& - (\Lambda^{1})^{T}L_{mnq}\Lambda^{1}  - (\Lambda^{2})^{T}L_{mnq}\Lambda^{2} \\
& - (x - x^{*})^{T}(f_{x} - f_{x^{*}})
\end{split}\label{Beta_E}
\end{gather}

Since $L_{mnq} \geq 0$ and $(x - x^{*})^{T}(f_{x} - f_{x^{*}}) \geq 0$ followed by the convexity of $f$, then according to Lemma 2.2, it follows from ($\ref{Beta_E}$) that (\ref{BA}) is satisfied.
\end{proof}

The following theorem proofs the convergence of trajectory $x(t)$ with the proposed algorithm ($\ref{Algorithm}$) to the optimal solutions.

\textbf{Theorem 5.1.} 
For Algorithm ($\ref{Algorithm}$) with Assumption 3.1, we have that the results that

(i) the trajectory $(x$, $Z$, $W$, $\Lambda^{1}$, $\Lambda^{2}$, $\phi)$ is bounded;

(ii) $x(t)$ converges to the optimal solution to Problem ($\ref{Problem}$).

\begin{proof}
i) Let $V(\phi)$ be as defined in ($\ref{Lyapunov}$). It follows from Lemma 5.1 that
\begin{align}
& \max \mathcal{L_{F}}V(\phi) \notag \\
\leq & \max \lbrace - (x - x^{*})^{T}(f_{x} - f_{x^{*}}) - (\Lambda^{1})^{T}L_{mnq}\Lambda^{1} \notag \\
& - (\Lambda^{2})^{T}L_{mnq}\Lambda^{2} : f_{x} \in \partial f(x), f_{x^{*}} \in \partial f(x^{*}), \label{LF_V}  \\
& x = P_{\Omega}[\bar{x}], \Lambda^{1} = P_{\bar{\mathbb{R}}_{+}^{mnq}}[\bar{\Lambda}^{1}], \Lambda^{2} = P_{\bar{\mathbb{R}}_{+}^{mnq}}[\bar{\Lambda}^{2}] \rbrace \leq 0 \notag
\end{align}

Note that $V(\phi) \geq \frac{1}{2}(\Vert x - x^{*} \Vert^{2} + \Vert Z - Z^{*} \Vert^{2} + \Vert W - W^{*} \Vert^{2} + \Vert U - U^{*} \Vert^{2} + \Vert \Lambda^{1} - \Lambda^{1*} \Vert^{2} + \Vert \Lambda^{2} - \Lambda^{2*} \Vert^{2} + \Vert Y^{1} - Y^{1*} \Vert^{2} + \Vert Y^{2} - Y^{2*} \Vert^{2})$
according to Lemma 2.1. Hence that trajectory ($x(t), Z(t), W(t), U(t), \Lambda^{1}(t), \Lambda^{2}(t),$ $Y^{1}(t), Y^{2}(t)$), $t \geq 0$ is bounded.

Because $\partial f(x)$ is compact for all $x \in \Omega$ and $(x(t), Z(t), W(t), U(t), \Lambda^{1}(t), \Lambda^{2}(t),Y^{1}(t), Y^{2}(t))$ is bounded for all $t \geq 0$, there exists $M = M(x, Z, W, \Lambda^{1}, \Lambda^{2}, \phi)> 0$ such that
\begin{gather}
\begin{cases}
M \geq & \!\!\!\! \Vert x(t) - f_{x(t)} - EA^{*}\Lambda^{1}(t) - E\hat{A}^{*}\Lambda^{2}(t) \Vert  \\
M \geq & \!\!\!\! \Vert Z(t) - \frac{1}{n}\Gamma \Lambda^{1}(t) + \Lambda^{2}(t) - L_{mnq} U(t) \Vert  \\
M \geq & \!\!\!\! \Vert W(t) - \Lambda^{1}(t) + \Lambda^{2}(t) \Vert  \\
M \geq & \!\!\!\! \Vert \Lambda^{1}(t) + [A^{*}E^{T}x(t) + \frac{1}{n}\Gamma Z(t) + W(t)\\
&  - B]  + L_{mnq}Y^{1}(t) - L_{mnq}\Lambda^{1} \Vert  \\
M \geq & \!\!\!\! \Vert \Lambda^{2}(t) + [\hat{A}^{*}E^{T}x(t) - Z(t) - W(t)]  \\
& + L_{mnq}Y^{2}(t) - L_{mnq}\Lambda^{2}(t) \Vert  
\end{cases}
\end{gather}
for all $f_{x(t)} \in \partial f(x(t))$ and all $t \geq 0$. Define $X: \mathcal{R}^{nq} \times \mathcal{R}^{mnq} \times \mathcal{R}^{mnq} \times \mathcal{R}^{mnq} \times \mathcal{R}^{mnq} \to \mathcal{R}$ by 
\begin{gather}
X\!(\!\bar{x}\!,\! \bar{Z}\!,\! \bar{W}\!,\! \bar{\Lambda}^{1}\!,\! \bar{\Lambda}^{2}\!) \!\!=\!\! \frac{1}{2}(\Vert \bar{x} \Vert^{2} \!\!\!+\!\! \Vert \bar{Z} \Vert^{2} \!\!\!+\!\! \Vert \bar{W} \Vert^{2} \!\!\!+\!\!\! \sum_{i=1}^{2}\Vert \bar{\Lambda}^{i} \Vert^{2})
\end{gather}

The function $X(\bar{x}, \bar{Z}, \bar{W}, \bar{\Lambda}^{1}, \bar{\Lambda}^{2})$ along the trajectories of ($\ref{Algorithm}$) satisfies that
\begin{gather}
\begin{split}
& \mathcal{L_{F}}X(\bar{x}, \bar{Z}, \bar{W}, \bar{\Lambda}^{1}, \bar{\Lambda}^{2})  \\
= & \lbrace x^{T}(- \bar{x} + x - f_{x} - EA^{*}\Lambda^{1} - E\hat{A}^{*}\Lambda^{2})   \\
&  + \bar{Z}^{T}(- \bar{Z} + Z - \frac{1}{n}\Gamma \Lambda^{1} + \Lambda^{2} - L_{mnq} U)   \\
& + \bar{W}^{T}(- \bar{W} + W - \Lambda^{1} + \Lambda^{2}) + (\bar{\Lambda}^{1})^{T}(\!-\! \bar{\Lambda}^{1} \!+\! \Lambda^{1} \\
& + [A^{*}\!E^{T}x \!+\! \frac{1}{n}\Gamma Z \!+\! W \!-\! B] \!+\! L_{mnq}Y^{1} \!\!\!-\! L_{mnq}\Lambda^{1}\!)   \\
& + (\bar{\Lambda}^{2})^{T}\!(\!-\! \bar{\Lambda}^{2} \!+\! \Lambda^{2} \!+\! [\hat{A}^{*}\!E^{T}x \!-\! Z \!-\! W] \!+\! L_{mnq}Y^{2} \\
& \!-\! L_{mnq}\Lambda^{2}\!) \!:\! 
f_{x} \!\in\! \partial f(x),\! x \!=\! P_{\Omega}[\bar{x}],\! Z \!=\! P_{\bar{\mathbb{R}}_{+}^{mnq}}[\bar{Z}],\\
&  W \!\!=\!\! P_{\bar{\mathbb{R}}_{+}^{mnq}}[\bar{W}], \Lambda^{1} \!\!=\!\! P_{\bar{\mathbb{R}}_{+}^{mnq}}[\bar{\Lambda}^{1}], \Lambda^{2} \!\!=\!\! P_{\bar{\mathbb{R}}_{+}^{mnq}}[\bar{\Lambda}^{2}] \rbrace  
\end{split}
\end{gather}

Note that
\begin{gather}
\begin{cases}
- \Vert \bar{x} \Vert^{2} \!\!+\! M\Vert \bar{x} \Vert \geq &\!\!\!\!\!\!\!\!\!\! x^{T}(t)(- \bar{x}(t) + x(t) - f_{x(t)}   \\
&\!\!\!\!\!\!\!\!\!\! - EA^{*}\Lambda^{1}(t) - E\hat{A}^{*}\Lambda^{2}(t))   \\
- \Vert \bar{Z} \Vert^{2} \!\!+\! M\Vert \bar{Z} \Vert \geq &\!\!\!\!\!\!\!\!\!\! \bar{Z}^{T}(t)(- \!\bar{Z}(t) \!\!+\!\! Z(t) \!\!-\!\! \frac{1}{n}\Gamma \Lambda^{1}(t)\\
&\!\!\!\!\!\!\!\!\!\! + \Lambda^{2}(t) - L_{mnq} U(t))   \\
- \Vert \bar{W} \Vert^{2} \!\!+\! M\Vert \bar{W} \Vert \geq &\!\!\!\!\!\!\! \bar{W}^{T}(t)(- \bar{W}(t) + W(t) \\
&\!\!\!\!\!\!\!\!\!\! - \Lambda^{1}(t) + \Lambda^{2}(t))   \\
- \Vert \bar{\Lambda}^{1} \Vert^{2} \!\!+\! M\Vert \bar{\Lambda}^{1} \Vert \geq &\!\!\!\!\!\! (\bar{\Lambda}^{1})^{T}(t)(-\! \bar{\Lambda}^{1}(t) \!+\! \Lambda^{1}(t)   \\
&\!\!\!\!\!\!\!\!\!\! + [A^{*}\!E^{T}x(t)\!\! +\! \frac{1}{n}\Gamma Z(t) \!+\! W(t)\\
&\!\!\!\!\!\!\!\!\!\! -\! B] \!\!+\!\! L_{mnq}Y^{1}\!(t) \!-\! L_{mnq}\Lambda^{1}\!(t))   \\
- \Vert \bar{\Lambda}^{2} \Vert^{2} \!\!+\! M\Vert \bar{\Lambda}^{2} \Vert \geq &\!\!\!\!\!\! (\bar{\Lambda}^{2})^{T}(t)(- \bar{\Lambda}^{2}(t) + \Lambda^{2}(t)   \\
&\!\!\!\!\!\!\!\!\!\! + [\hat{A}^{*}E^{T}x(t) - Z(t) - W(t)]  \\
&\!\!\!\!\!\!\!\!\!\! + L_{mnq}Y^{2}(t) - L_{mnq}\Lambda^{2}(t)) 
\end{cases}
\end{gather}

Hence,
\begin{gather}
\begin{split}
& \max \mathcal{L_{F}}X(\bar{x}(t), \bar{Z}(t), \bar{W}(t), \bar{\Lambda}^{1}(t), \bar{\Lambda}^{2}(t))   \\
\leq & -2X(\bar{x}(t), \bar{Z}(t), \bar{W}(t), \bar{\Lambda}^{1}(t), \bar{\Lambda}^{2}(t))   \\
& + 5M \sqrt{X(\bar{x}(t), \bar{Z}(t), \bar{W}(t), \bar{\Lambda}^{1}(t), \bar{\Lambda}^{2}(t))} 
\end{split}
\end{gather}

It can be easily verified that $X(\bar{x}(t)$, $ \bar{Z}(t)$, $\bar{W}(t)$, $\bar{\Lambda}^{1}(t)$, $\bar{\Lambda}^{2}(t))$, $t \geq 0$, is bounded, so are $\bar{x}(t)$, $\bar{Z}(t)$, $\bar{W}(t)$, $\bar{\Lambda}^{1}(t)$, $\bar{\Lambda}^{2}(t)$ for all $t \geq 0$. As the result, the trajectory $(x$, $Z$, $W$, $\Lambda^{1}$, $\Lambda^{2}$, $\phi)$ is bounded.

ii) Let 
\begin{gather}
\begin{split}
\mathcal{R} \subset & \lbrace \phi \in \mathbb{R}_{(7m+1)nq} :       x = P_{\Omega}[\bar{x}], x^{*} = P_{\Omega}[\bar{x}^{*}], \\
& \min_{f_{x} \in \partial f(x), f_{x^{*}} \in \partial f(x^{*})} (x - x^{*})^{T}(f_{x} - f_{x^{*}}) = 0,   \\
& L_{mnq}\Lambda^{1} = \textbf{0}_{mnq}, L_{mnq}\Lambda^{2} = \textbf{0}_{mnq} \rbrace  
\end{split}
\end{gather}

Note that $(x - x^{*})^{T}(f_{x} - f_{x^{*}}) > 0$ if $x \neq x^{*}$ since the Assumption 3.1. Hence, $\mathcal{R} \subset \lbrace \phi \in \mathbb{R}_{(7m+1)nq} : L_{mnq}\Lambda^{1} = \textbf{0}_{mnq}, L_{mnq}\Lambda^{2} = \textbf{0}_{mnq}, x = P_{\Omega}[\bar{x}] = x^{*} \rbrace$. Let $\mathcal{M}$ be the largest weakly invariant subset of $\bar{\mathcal{R}}$. According to Lemma 2.3, $\phi \to \mathcal{M}$ as $t \to \infty$. Hence, $x(t) \to x^{*}$ as $t \to \infty$. Part (ii) is thus proved.
\end{proof}

\section{Simulation}

In this section, we show a numerical example to validate our proposed distributed optimization algorithm. Consider the distributed robust optimization problem with four agents moving in a 2-D space with first-order dynamics as follows
\begin{gather}
F(x)= \sum_{i=1}^{4} \Vert x_{i}-p_{i} \Vert^{2}_{2} + \vert x \vert_{1}
\end{gather}
where $p_{i}=[i,-i]^{T}$, $\Vert \cdot \Vert_{2}$ denotes the $l_{2}$ norm, $\Omega_{i} = \lbrace \delta \in \mathbb{R}^{2} \vert \Vert \delta - x_{i}(0) \Vert_{2} \leq 30\rbrace$, $m=2$, $\gamma_{1}=\gamma_{2}=2$, $A_{i1}=0.1 \cdot i \cdot I_{2}$, $\hat{A}_{i1}=0.1 \cdot (5-i)\cdot I_{2}$, $A_{i2}=\hat{A}_{i1}$, $\hat{A}_{i2}=A_{i1}$, and $b_{1}^{1}=[-15,-5]^{T}$, $b_{1}^{2}=[-10,-4]^{T}$, $b_{1}^{3}=[0,-6]^{T}$, $b_{1}^{4}=[4,0]^{T}$, $b_{2}^{1}=[-5,-1]^{T}$, $b_{2}^{2}=[-4,-3]^{T}$, $b_{2}^{3}=[0,-2]^{T}$, $b_{2}^{4}=[1,-5]^{T}$, $b_{1}=[-21,-15]^{T}$, $b_{2}=[-8,-11]^{T}$. This problem can be transferred to its corresponding dual problem as the form of problem (\ref{Dual-Problem}). The Laplacian of the undirected graph $\mathcal{G}$ is given by
\begin{equation}
L_{4} = 
\begin{bmatrix}
 1   & -1  &  0  & 0   \\
-1   &  2  & -1  & 0   \\
 0   & -1  &  2  & -1  \\
 0   &  0  & -1  & 1   \\
\end{bmatrix}
\end{equation} 
The initial positions of the agents 1, 2, 3, and 4 are set as $x_{1}(0) = [-13, 12]^{T}$, $x_{2}(0) = [17, 15]^{T}$, $x_{3}(0) = [-10, -11]^{T}$ and $x_{4}(0) = [16, -14]^{T}$. We set the initial values for the Lagrangian multipliers $\lambda_{ij}^{1}$, $\lambda_{ij}^{2}$, $\mu_{ij}$  and auxiliary variables $z_{ij}$, $w_{ij}$, $y_{ij}^{1}$, $y_{ij}^{2}$ as zeros for $i \in \lbrace 1, 2, 3, 4 \rbrace, j \in \lbrace 1, 2 \rbrace$. The optimal solution is $x^{*}_{1}=[-7.439, -10.408]^{T}, x^{*}_{2}= [-4.016, -6.409]^{T}, x^{*}_{3}=[-15.516, -17.612]^{T}, x^{*}_{4}=[-13.401, -19.965]^{T}$.

\begin{figure}
\centering
\subfigure{
\includegraphics[width=0.2\textwidth]{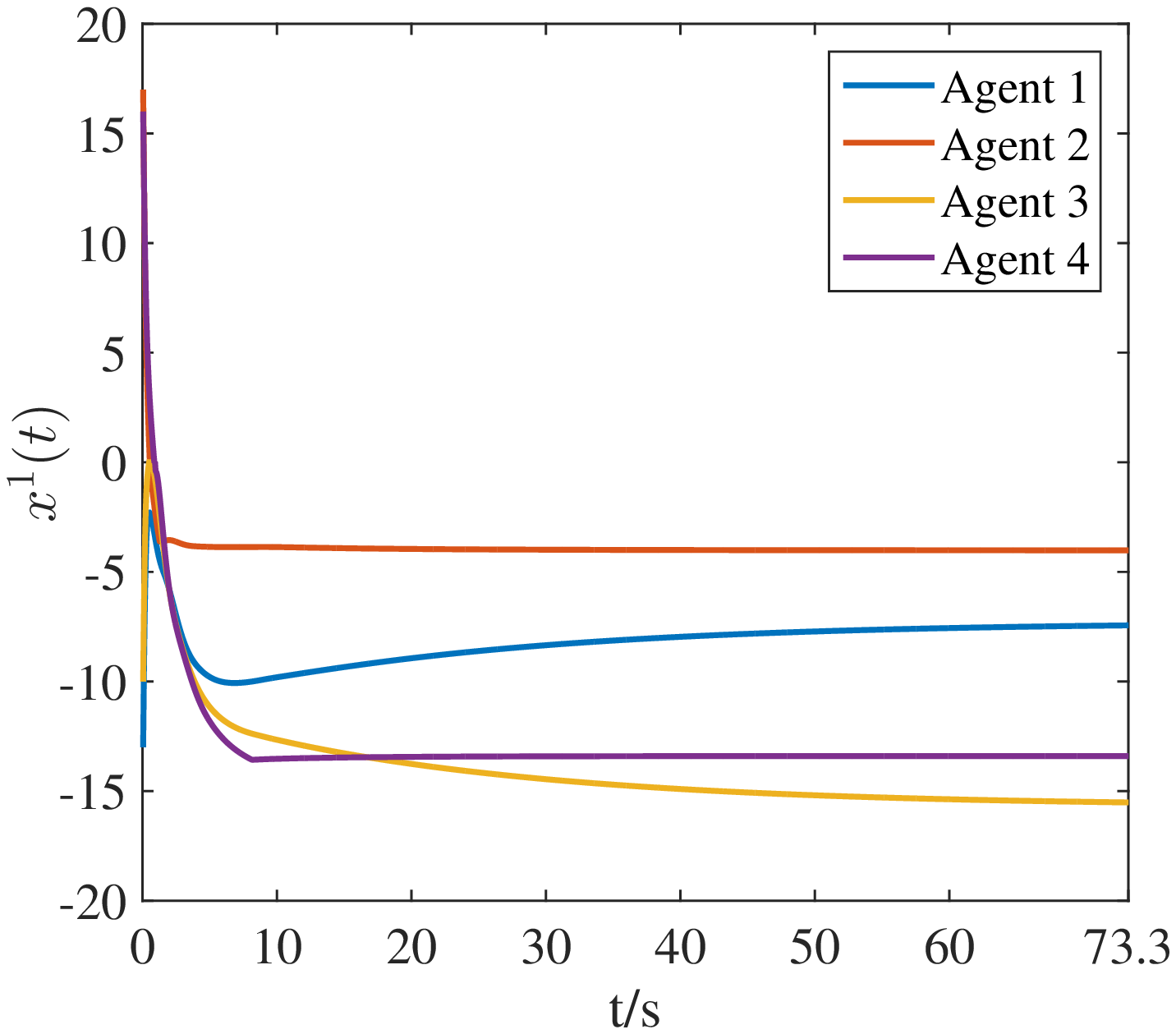}}
\subfigure{
\includegraphics[width=0.2\textwidth]{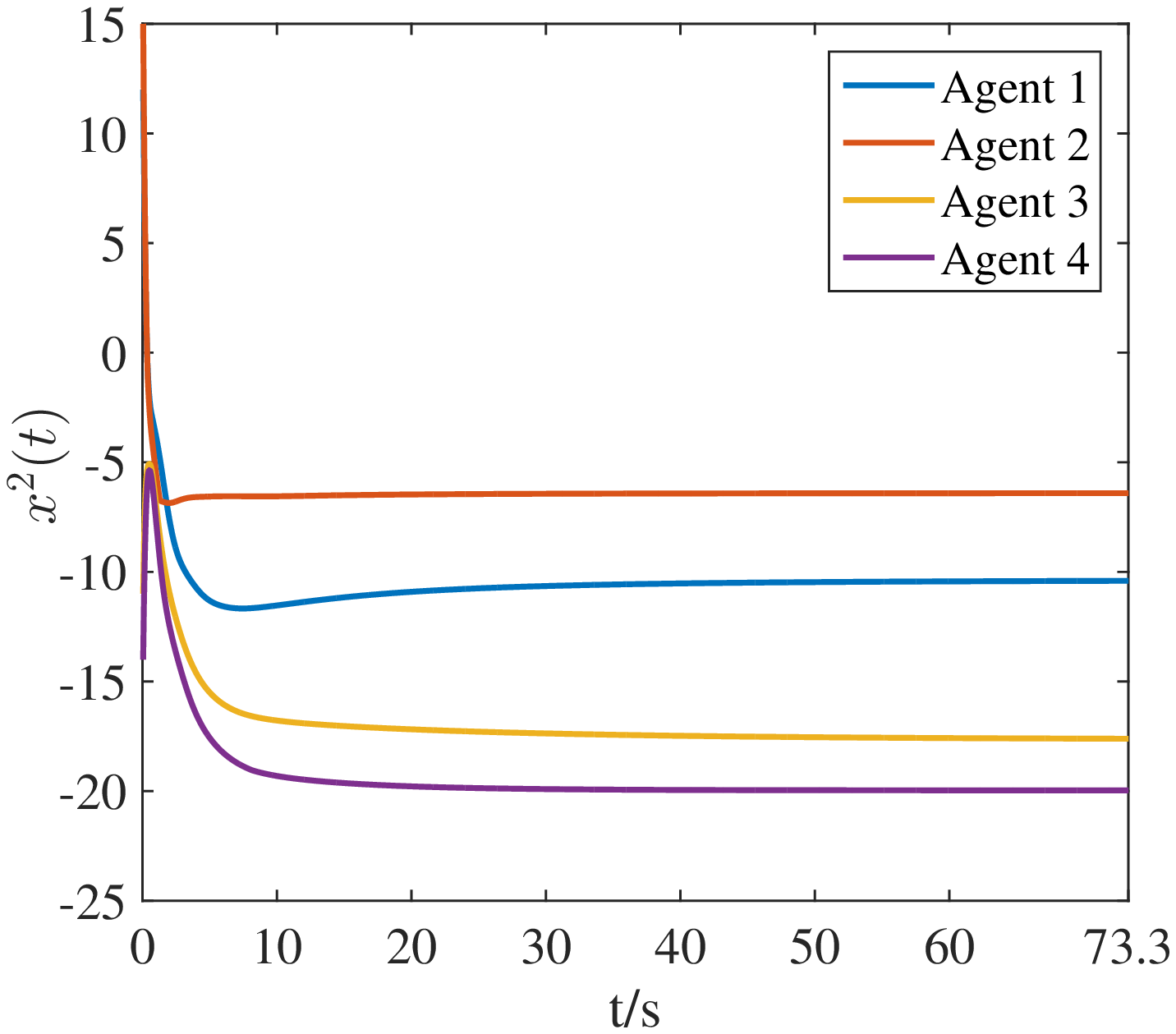}}
\caption{The trajectories of $x_{i}(t), i \in \lbrace 1, 2, 3, 4 \rbrace$ with algorithm ($\ref{Algorithm}$)}
\label{Fig.1}
\end{figure}

\begin{figure}
\centering
\subfigure{
\includegraphics[width=0.20\textwidth]{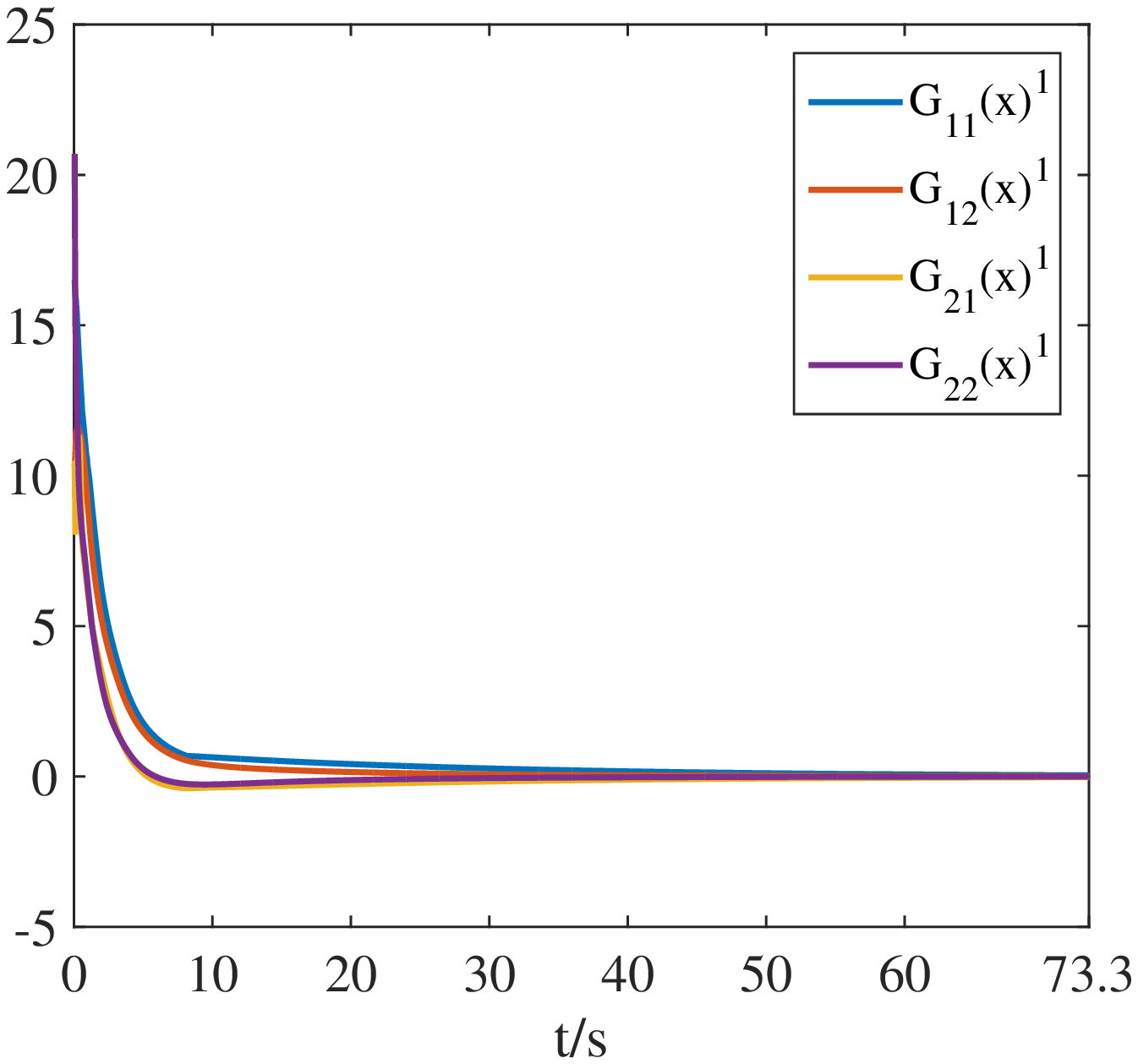}}
\subfigure{
\includegraphics[width=0.20\textwidth]{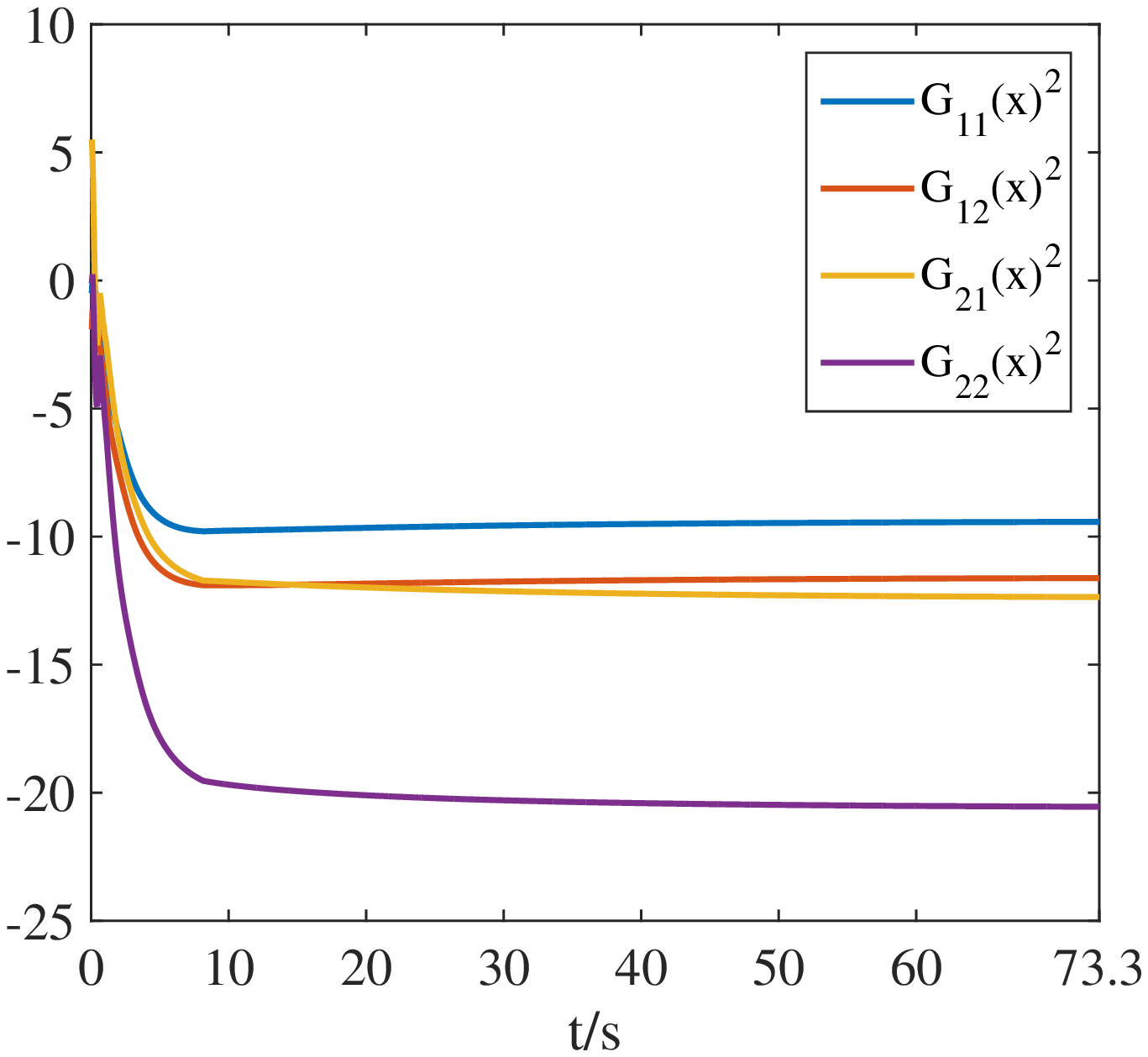}}
\caption{The trajectories of $G_{j1}(x)$ and $G_{j2}(x)$, $j \in \lbrace 1, 2\rbrace$ with algorithm ($\ref{Algorithm}$)}
\label{Fig.2}
\end{figure}

Fig.$\ref{Fig.1}$ gives the trajectories of $x_{i}(t), i \in \lbrace 1, 2, 3, 4 \rbrace$. It can be seen that the trajectory of $x$ converges to the optimal solution. Let $G_{j1}(x)=\sum_{i=1}^{4} H_{ij}^{1}$, $G_{j2}=\sum_{i=1}^{4} H_{ij}^{2}$, $j \in \lbrace 1, 2\rbrace$.  Fig.$\ref{Fig.2}$ shows the trajectory of $G_{j1}(x)$ and $G_{j2}(x)$, $j \in \lbrace 1, 2\rbrace$, which proves that the constraint condition of problem (\ref{Problem}) are satisfied.

\section{Conclusion}

In this paper, a distributed nonsmooth resource allocation problem with cardinality constrained uncertainty has been investigated. With the help of duality theory about convex optimization, a deterministic distributed robust resource allocation problem with linear optimization formulation has been derived under the framework of multi-agent system. A distributed projection-based algorithm has been proposed to deal with this problem. Based on stability theory and differential inclusions, the proposed algorithm has been proved to reach the optimal solution and satisfy the resource allocation condition 
simultaneously.

\end{document}